\documentclass[preprint,12pt]{elsarticle}
\newtheorem{definition}{\textbf{Definition}}

\newtheorem{theorem}{\textbf{Theorem}}

\def\Q {\mathbb{Q}}
\def\R {\mathbb{R}}

\usepackage{amsfonts}
\usepackage{amssymb}
\def\QQ {\overline{\Q}}

\usepackage{mathrsfs}

\def\S {\mathcal{S}}
\def\cH {\mathcal{H}}

\def\Z {\mathbb{Z}}

\def\C {\mathbb{C}}

\def\dfrac {\displaystyle\frac}

\def\den{\mathop{\rm den}}

\def\tareesidedbox#1{\setbox0=\hbox{$#1$}\dimen0=\wd0 \advance\dimen0 by3pt\rlap{\hbox{\vrule height8pt width.4pt
 depth2pt \kern-.4pt\vrule height8.4pt width\dimen0 depth-8pt\kern-.4pt \vrule height8pt width.4pt depth2pt}}
 \relax \hbox to\dimen0{\hss$#1$\hss}}

\usepackage{amssymb}

\journal{...}

\begin{document}

\begin{frontmatter}



\title{On transcendental analytic functions mapping an uncountable class of $U$-numbers into Liouville numbers}



\author{Diego Marques\fnref{dm}\corref{CA}}{\ead{diego@mat.unb.br}}
\author{Josimar Ramirez\fnref{jr}}{\ead{jrmat0702@gmail.com}}
\address{Departamento de Matem\' atica, Universidade  de Bras\' ilia, Bras\' ilia, 70910-900, Brazil}


\fntext[dm]{Supported by FAP-DF and CNPq-Brazil}
\fntext[jr]{Scholarship holder of CAPES-Brazil}

\cortext[CA]{Corresponding author}

\begin{abstract}
In this paper, we shall prove, for any $m\geq 1$, the existence of an uncountable subset of $U$-numbers of type $\leq m$ (which we called the set of {\it $m$-ultra numbers}) for which there exists uncountably many transcendental analytic functions mapping it into Liouville numbers.
\end{abstract}


\end{frontmatter}

\section{Introduction}

A \textit{transcendental function} is a function $f(x)$ such that the only complex polynomial satisfying $P(x, f(x)) =0$, for all $x$ in its domain, is the null polynomial. For instance, the trigonometric functions, the exponential function, and their inverses. 

The study of the arithmetic behavior of transcendental functions at complex points has attracted the attention of many mathematicians for decades. The first result concerning this subject goes back to 1884, when Lindemann proved that the transcendental function $e^z$ assumes transcendental values at all nonzero algebraic point. In 1886, Weierstrass gave an example of
a transcendental entire function which takes rational values at all rational points. Later, St$\ddot{\mbox{a}}$ckel \cite{19} proved
that for each countable subset $\Sigma\subseteq \C$ and each dense subset $T\subseteq \C$, there is a transcendental entire
function $f$ such that $f(\Sigma) \subseteq T$ (F. Gramain showed that St$\ddot{\mbox{a}}$ckel's theorem is valid if $\Sigma$ and $T$ are subsets of $\R$). Another construction due to St$\ddot{\mbox{a}}$ckel \cite{20} produces a transcendental entire function $f$ whose
derivatives $f^{(t)}$, for $t = 0, 1, 2, \ldots$, all map $\Q$ into $\Q$. Two years later, G. Faber refined this result by showing the existence of a transcendental entire function such that $f ^{(t)}(\QQ)\subseteq  \Q(i)$, for all $t \geq 0$. In 1968, van der Poorten \cite{alf} proved the existence of a transcendental function $f$, such that $f^{(s)}(\alpha)\in \Q(\alpha)$, for all $\alpha\in \QQ$. In 2011, Marques \cite{Diego} proved some of these results in the hypertranscendental context (for more on this subject, we refer the reader to \cite{wal1} and the references therein). Let $\den(z)$ be the denominator of the rational number $z$. Very recently, Marques and Moreira \cite{gugu} proved the existence of uncountable many transcendental entire functions $f$, such that $f(\Q)\subseteq \Q$ and $\den(f(p/q))<q^{8q^2},$ for all $p/q\in \mathbb{Q}$, with $q>1$. 

A real number $\xi$ is called a \textit{Liouville number}, if there exists a rational sequence $(p_k/q_k)_{k\geq 1}$, with $q_k> 1$, such that
\begin{center}
$0<\left|\xi - \dfrac{p_k}{q_k}\right|<q_k^{-k}$, for $k=1,2,\ldots$.
\end{center}
The set of the Liouville numbers is denoted by $\mathbb{L}$ and it is a dense $G_{\delta}$ set (and therefore uncountable).

All the previously mentioned results deal with the arithmetic behavior of a countable set by a transcendental function. However, in 1984, in one of his last papers, K. Mahler \cite{Mahler} raised the following question:

\medskip

\noindent
{\bf Question.} Are there transcendental entire functions $f(z)$ such that if $\xi$ is any Liouville number, then so is $f(\xi)$?
\medskip

He also said that: ``The difficulty of this problem lies of course in the fact that the set of all Liouville numbers is non-enumerable". In this direction, Marques and Moreira \cite{gugu} showed the existence of an uncountable subset of Liouville numbers for which there exists uncountably many transcendental entire functions mapping it into itself.

As usual, $\cH(\alpha)$ will denote the height of the algebraic number $\alpha$ (that is, the maximum of the absolute values of its primitive minimal polynomial over $\Z$) and $\exp^{[3]}(x)=e^{e^{e^x}}$. Now, let us define the following class of numbers:

\begin{definition}
A real number $\xi$ is called an \emph{$m$-ultra number} if there exist infinitely many $m$-degree real algebraic numbers $(\alpha_n)_n$ such that
\begin{center}
$0<\left|\xi - \alpha_n\right|<(\exp^{[3]}(\cH(\alpha_n)))^{-n}$, for $n=1,2,\ldots$.
\end{center}
The set of the $m$-ultra numbers will be denoted by $\mathbb{U}_{\scriptsize{\emph{$m$-ultra}}}$.
\end{definition}

It follows from the definition that $\mathbb{U}_{\scriptsize{\mbox{$m$-ultra}}}$ is a subset of $U$-numbers of type at most $m$ and it is also a dense $G_{\delta}$ set (in particular it is uncountable) - it means that $\mathbb{U}_{\scriptsize{\emph{$m$-ultra}}}$ is a large set in a topological sense. In particular, every real number can be written as the sum (and the product, provided it is non zero) of two $m$-ultra numbers, as in \cite{erdos}. 

The aim of this paper is to prove the following result:

\begin{theorem}\label{1}
There exist uncountable many analytic transcendental functions $\phi:\R\to \R$ such that $\phi(\mathbb{U}_{\scriptsize{\emph{$m$-ultra}}})\subseteq \mathbb{L}$.
\end{theorem}

Recall that $\QQ_m$ denotes the set of all $m$-degree real algebraic numbers. Since $\QQ_m$ and $\Q$ are dense countable sets of $\R$, there exist uncountable many transcendental analytic functions $\phi$ with $\phi(\overline{\Q}_m) \subseteq \Q$. In order to prove Theorem \ref{1} we shall find a class of such functions with an upper bound for $\den (\phi(\alpha))$ in terms of $m$ and $\cH(\alpha)$. More precisely, we have

\begin{theorem}\label{2}
For any given positive integer $m\geq 1$, there exist uncountably many transcendental analytic functions $\phi:\R \to \R$ with $|\phi'(x)|<0.0001$, $\phi(\overline{\Q}_m) \subseteq \Q$ and such that for all $ \alpha \in \overline{\Q}_m$, it holds that
\begin{equation}\label{<q}
\den (\phi(\alpha)) \leq  (2q)^{450m^52^{18m^2}q^{6m}},
\end{equation}
where $q=\cH(\alpha)$. 
\end{theorem}

\section{The proofs}\label{sec2}
\subsection{Proof that Theorem \ref{2} implies Theorem \ref{1}}

Given an $m$-ultra number $\xi$, there exist infinitely many $\alpha_n\in \QQ_m$, with height at least $\max\{m,8\}$, and such that
\begin{center}
$0<\left|\xi - \alpha_n\right|<\dfrac{1}{(\exp^{[3]}(\cH(\alpha_n))^n}$, for $n=1,2,\ldots$.
\end{center}

Let $\phi$ be a function as in Theorem \ref{2}. By the Mean Value Theorem, we obtain
\[
\left|\phi(\xi) - \phi(\alpha_n)\right|\leq 0.0001\cdot \left|\xi - \alpha_n\right|<\dfrac{1}{(\exp^{[3]}(H(\alpha_n)))^n}.
\]
We know that $\phi(\alpha_n)=p_n/q_n$, with $q_n\leq (2t_n)^{450m^52^{18m^2}t_n^{6m}}$, where $t_n=\cH(\alpha_n)$. Since $t_n\geq \max\{m,8\}$, then a straightforward calculation gives $q_n\leq \exp^{[3]}(t_n)$ and hence
\begin{center}
$\left|\phi(\xi) - \dfrac{p_n}{q_n}\right|=\left|\phi(\xi) - \phi(\alpha_n)\right|<\dfrac{1}{q_n^n}$, for $n=1,2,\ldots$.
\end{center}
This implies that $\phi(\xi)$ is a Liouville number as desired.\qed

\subsection{Proof of Theorem \ref{2}}

Before starting the proof, we shall state some useful facts 
\begin{itemize}
\item For any distinct $y,b\in [-1,1]$, it holds that $|\sin (y-b)|>|y-b|/3$. \hfill \break
(Indeed, the function $\sin(x)/x$ is decreasing for $x\in (0,\pi]$, and $\sin(2)/2>1/3$.)

\item For any $x,y \in \overline{\Q}_m$ we have $\cH(y-x) \leq 2^{4m^2} \cH(x)^m \cH(y)^m$ \hfill \break
(Indeed, let $W(x)$ be the absolute multiplicative Weil height, we know that 
\begin{equation}\label{W}
	\frac{1}{2^{\partial (x)}} W(x)^{\partial(x)} \leqslant \cH(x) \leqslant 2^{\partial(x)} W(x)^{\partial(x)},
\end{equation}
where $\partial(x)$ is the degree of the algebraic number $x$. Using this together with the inequality $W(x+y) \leq
2 W(x)W(y)$ the result follows.)

\item For any distinct $x,y\in \overline{\Q}_m \cap [0,1/2]$, with $\cH(x),\cH(y)\leq n$, we have 
\[
|\cos (\pi x)-\cos (\pi y)|\geq \frac{\pi}{2^{4m^2+1} n^{2m+1}}.
\]
(Indeed, we can assume $x<y$. Then by the mean value theorem, one has
$
	|\cos (\pi x)-\cos (\pi y)|
	\geq
	\sin(\pi x) (\pi y- \pi x)
	\geq
	2\pi x(y-x).
$ 
There is a simple lower bound for the modulus of a nonzero complex algebraic numbers $\alpha$ in terms of the 
height $\cH(\alpha)$, namely $|\alpha| \geq (\cH(\alpha)+1)^{-1}\geq (2\cH(\alpha))^{-1}$ (see \cite[Page 82]{wal2}). Thus,
\begin{eqnarray*}
	|\cos (\pi x)-\cos (\pi y)| & > & \frac{\pi}{2H(x)H(y-x)}\\
	& \geq & 
	\frac{\pi}{2^{4m^2+1} n^{2m+1}}.)
\end{eqnarray*}

\item For every $\epsilon>0$, any interval of length $>\epsilon$ contains at least two rational numbers with denominator $\leq \lceil 2/\epsilon \rceil$.\hfill\break
(Indeed, if $m=\lceil 2/\epsilon \rceil$ and $(a,b)$ is the interior of the interval, we have $b-a>\epsilon\ge 2/m$, and so, for $k=\lfloor ma \rfloor+1$, we have $ma<k\le ma+1$, and so $ma<k<k+1\le ma+2<ma+m(b-a)=mb$, which implies $a<k/m<(k+1)/m<b$.)

\item Let $f:\mathbb{R}\to \mathbb{R}$ be a periodic function which assumes infinitely many values. Then $f$ is transcendental.\hfill\break
(Indeed, suppose that $f$ is algebraic, and let $P(x,y)=\sum_{i=0}^na_i(x)y^i$ be the polynomial with minimal degree $n$ in the variable $y$, such that $a_n(x)$ has minimal degree and $P(x,f(x))=0$, $\forall x\in \R$. We may suppose that $a_n(x)$ is monic. Since $f$ assumes infinitely many values, then $a_0(x),\ldots, a_n(x)$ cannot be all constants. Set $\ell$ as the largest index with $a_{\ell}(x)$ non constant. If $t$ is the period of $f$, then $Q_k(x,f(x)):=P(x+tk,f(x+tk))-P(x,f(x))=\sum_{i=0}^{\ell}(a_i(x+tk)-a_i(x))(f(x))^i=0$, for all $(x,k)\in \R\times \Z$. Note that for some integer $k_0$, $a_{\ell}(x+tk_0)-a_{\ell}(x)$ is nonzero. If $\ell=n$, then $a_n(x+tk)-a_n(x)$ is nonzero having degree smaller than the degree of $a_n(x)$ which contradicts the minimality of the degree of $a_n(x)$. In the case of $\ell<n$, then $Q_{k_0}(x,f(x))=0$ and $Q_{k_0}$ has degree $\ell<n$, in $y$, which contradicts our assumption on the minimality of $n$.)


\item Let $f:\R\to \R$ be a transcendental function and let $g:\R\to \R$ be a non constant algebraic function. Then $f\circ g$ is transcendental.\hfill\break
(Indeed, by assumption, $\C(y,f(y))$ is transcendental over $\C(y)$. Setting $y=g(x)$ is merely making an algebraic extension of each, so $\C(g(x),f(g(x)))$ is transcendental over $\C(g(x))$. Thus $\C(g(x),f(g(x)))$ is transcendental over $\C(g(x))$. Hence the tower $\C(g(x),f(g(x)))\supset \C(g(x))\supset \C(x)$ is transcendental, so $f(g(x))$ is transcendental over $\C(x)$.)

\end{itemize}


Now, we are ready to deal with the proof of the theorem.

Consider the following enumeration of $A:=\overline{\Q}_m \cap [0,1/2]$:
$$
	A=\{\alpha_1,\alpha_2,\alpha_3,\ldots \},
$$
constructed as follows. Let $\S_k$ be the set of all irreducible and primitive polynomials in $\Z[x]$ with degree $m$ and height $k$. Denote by $t_k:=|\S_k|<(m+1)(2k+1)^m$. Let $\mathcal{R}_k$ be the set all distinct roots of polynomial in $\S_k$ belonging to the interval $[0,1/2]$, (Note that $\mathcal{R}_k \cap \mathcal{R}_t = \emptyset$, for $k \neq t$) and 
$l_k=|\mathcal{R}_k|$, then $\mathcal{R}_k =\{\gamma_1^{(k)}, \ldots, \gamma_{l_k}^{(k)} \}$ with 
$\gamma_i^{(k)} < \gamma_{i+1}^{(k)}$ $\forall k \geq 1$. So the desired enumeration is given by
$$
	A=\{\alpha_1,\alpha_2,\alpha_3,\ldots \}= \{ \mathcal{R}_1, \mathcal{R}_2, \mathcal{R}_3, \ldots \}.
$$
Now, we will give estimates for the height of the algebraic numbers in $A$ as a function of the position in the enumeration. Although estimates are not the best, they will be sufficient for our purposes.

If $\alpha_n \in \mathcal{R}_{k+1}$, then $\cH(\alpha_n)= k+1$. We have $n\leq l_1+\cdots+l_{k+1} \leq
t_1+ \cdots + t_{k+1} \leq (m+1)(2k+3)^{m+1}$, therefore
$$
	\cH(\alpha_n) \geqslant \frac{1}{2} \sqrt[m+1]{\frac{n}{m+1}} -2.
$$
On another hand, $n\geq l_1+\cdots+ l_k$. Let $j$ be an odd number with $4<j \leq k$, then $l_j\geq 1$
(because $(2/j)^{1/m} \in \mathcal{R}_j$). Thus, if $k\geq 5$, we have $n \geq \lfloor \frac{k-4}{2}\rfloor > \frac{k-6}{2}$, therefore
$\cH(\alpha_n) < 2n+7$ (the cases $k=1,2,3,4$ are trivial). Define $B_n= \{y_1,y_2,\ldots,y_n\}$ with $y_k:= \cos(\pi \alpha_k)$.

Set $h:\C\to \C$ given by
\[
h(x)=g(\cos (\pi x)),
\]
where $g(y)=\sum_{n=1}^{\infty}c_ng_n(y)$, with $g_n(y)=\prod_{b\in B_n}\sin (y-b)$. 


Suppose that $c_n=0$ for $1\le n\le 5$ and $|c_n|<1/n^n$ for every positive integer $n$. We claim that $h$ is an entire function. In fact, for all $y$ belonging to the open ball $B(0,R)$ one has that
\[
|g_n(y)|<\displaystyle\prod_{b\in B_n}e^{|y-b|}\leq e^{n(R+1)},
\]
where we used that $b\in [-1,1]$. Thus, since $|c_n|<1/n^n$, we get $|c_ng_n(y)|\leq (e^{R+1}/n)^n$ yielding that $g$ (and so $h$) is an entire function, since the series $g(y)=\sum_{n=1}^{\infty}c_ng_n(y)$, which defines $g$, converges uniformly in any of these balls. Let $f:\R\to \R$ be the restriction of $h$ to $\R$. In particular, $f$ is analytic and $|f'(x)|\leq \sum_{n=6}^{\infty}1/n^{n-1}<0.0002$, for all $x\in \R$.

Now, we shall choose inductively  $c_n$'s conveniently such that $f$ satisfies $f(\alpha_k)\in  \Q$, for all $k$, and $\den (f(\alpha_k))< (72m^2(6q)^{4m})^{10m^3(6q)^{2m}}$, where $q=\mathcal{H}(\alpha_k)$. 

Suppose that $c_1,\ldots, c_{n-1}$ were chosen such that $f(\alpha_1),\ldots, f(\alpha_n)$ have the desired properties (notice that the choice of $c_1,\ldots, c_{n-1}$ determines the values of $f(\alpha_1),\ldots, f(\alpha_n)$, independently of the values of $c_k,k\ge n$; in particular, since $c_k=0$ for $1\le k\le 5$, we have $f(\alpha_n)=0$ for $1\le n\le 6$). Now, we shall choose $c_n$ for which $f(\alpha_{n+1})$ satisfies the requirements.

Let $t\leq n$ be positive integers with $n\ge 5$. Then $\cH(\alpha_{n+1}),\cH(\alpha_t)\leq 2n+9$. Since $\cos (\pi \alpha_{n+1})\neq \cos (\pi \alpha_{t})$, then 
$$
	|y_{n+1}-y_t|
	\geq
	\frac{\pi}{2^{4m^2+1}(2n+9)^{2m+1}}.
$$
Therefore
\[
|\sin (y_{n+1}-y_t)|>\frac{|y_{n+1}-y_t|}{3}>\frac{\pi/3}{2^{4m^2+1}(2n+9)^{2m+1}}
\]
yielding $|g_n(y_{n+1})|>(\frac{\pi/3}{2^{4m^2+1}(2n+9)^{2m+1}})^{n}$. Thus $c_ng_n(y_{n+1})$ runs an interval of length larger than $2\pi^n (3n)^{-n} 2^{-n(4m^2+1)} (2n+9)^{-3mn}$. Now, we may choose 
(in at least two ways) $c_n\neq 0$ such that $g(y_{n+1})$ is a rational number with denominator at most 
$n^n 2^{n(4m^2+1)}(2n+9)^{3mn}$. Thus $\den (f(\alpha_k))=\den(g(\cos(\pi \alpha_k))=\den(g(y_k))\leq (k-1)^{(k-1)} 2^{(k-1)(4m^2+1)} (2k+7)^{3m(k-1)} <k^k 2^{k(4m^2+1)} (2k+7)^{3mk}$. Since $q:=\cH(\alpha_k)\geq \frac{1}{2} \sqrt[m+1]{k/(m+1)}-2$, we get $k \leq (2q+4)^{(m+1)}(m+1)$. Then
\begin{eqnarray*}
 \den (f(\alpha_k))  & \geq & \left((2q+4)^{(m+1)}(m+1) \right)^{(2q+4)^{(m+1)}(m+1)} 2^{(2q+4)^{(m+1)}(m+1)(4m^2+1)} \\
&  & \times ( 2(2q+4)^{(m+1)}(m+1)+ 7 )^{3m(2q+4)^{(m+1)}(m+1)}\\
& < & (72m^2(6q)^{4m})^{10m^3(6q)^{2m}}.
\end{eqnarray*}

Now, consider the function $\psi:\R\to \R$, given by $\psi(x):=\frac{x}{2(1+x^2)}$. Note that $\psi(\QQ_m)\subseteq \QQ_m\cap [0,1/2]$. Therefore, our desired function is $\phi:=f\circ \psi$. In fact, $\phi(\QQ_m)\subseteq \Q$, $|\phi'(x)|=|f'(\psi(x))||\psi'(x)|< 0.0001$, for all $x\in \R$, and by our previous argument, if $\alpha\in \QQ_m$, then
\[
\den (\phi(\alpha)) =\den (f(\psi(\alpha)))\leq (72m^2(6t)^{4m})^{10m^3(6t)^{2m}},
\]
where $t=\mathcal{H}(\psi(\alpha))$. On the other hand, we can use (\ref{W}) together with the fact that $W(x/y)\leq W(x)W(y)$ to obtain $t=H(\frac{\alpha}{2(1+\alpha^2)})\leq 2^{6m}q^{3}$, where $q=\mathcal{H}(\alpha)$. Thus
\[
\den (\phi(\alpha)) < (2q)^{450m^52^{18m^2}q^{6m}}, 
\]
as desired.

Note that since there is a binary tree of different possibilities for $f$ (if we have choosen $c_1, c_2, \dots, c_{n-1}$, different choices of $c_n$ give different values of $f(y_{n+1})$, which does not depend on the values of $c_k$ for $k>n$, and so different functions $f$), we constructed uncountably many possible functions $f$. So there exist uncountable many functions $\phi$ (since $\psi$ is non constant).

Now, it remains to prove that all functions constructed above are transcendental: in fact, since $f$ assumes infinitely many values (because it is continue and non constant) and it is periodic (with period $2$), then $f$ is transcendental. Therefore $f\circ \psi$ is transcendental, because $\psi$ is a non constant rational function.
\qed

\section*{acknowledgements}
The authors would like to thank Joseph Silverman and Carlos Gustavo Moreira  for helpful suggestions which improved the quality of this paper. The first author is supported by FAP-DF and CNPq-Brazil and the second author is scholarship holder of CAPES-Brazil.



\end{document}